\documentclass{amsart}
\usepackage{amssymb}
\usepackage{amsmath}
\usepackage{amsfonts}
\usepackage{amsmath}

\usepackage{mathptmx}
\usepackage{helvet}
\usepackage{courier}
\usepackage{type1cm}
\usepackage{graphicx}
\usepackage{multicol}
\usepackage[bottom]{footmisc}
\usepackage{amsmath}
\usepackage{amssymb}
\usepackage{amscd}
\usepackage{indentfirst}
\usepackage{graphicx}
\usepackage{epsfig}
\usepackage{epstopdf}
\usepackage{algorithm}
\usepackage{algorithmicx}
\usepackage{algpseudocode}
\usepackage[table]{xcolor}
\usepackage{eurosym}
\usepackage{ragged2e}
\usepackage{float}
\usepackage{graphics}
\usepackage[font=bf,labelsep=period]{caption}
\usepackage{epstopdf}
\usepackage{hyperref}
\usepackage{subcaption}

\usepackage{listings,xcolor}
\newtheorem{theorem}{Theorem}
\theoremstyle{plain}

\newtheorem{corollary}{Corollary}

\newtheorem{remark}{Remark}

\numberwithin{equation}{section}

\begin{document}

\begin{center}
\thispagestyle{empty} \pagestyle{myheadings}

\title[Implementation of computation formulas for certain classes of Apostol-type poly...]{\Large Implementation of computation formulas for certain classes of Apostol-type polynomials and some properties associated with these polynomials}

\maketitle

\author{\large Irem Kucukoglu} \\[0pt]
\textit{Department of Engineering Fundamental Sciences, Faculty of
Engineering,\\[0pt]
Alanya Alaaddin Keykubat University TR-07425 Antalya, Turkey, \\[0pt]
}\textbf{\texttt{irem.kucukoglu@alanya.edu.tr}}

\bigskip
\end{center}

\textbf{Abstract}. The main purpose of this paper is to present various identities and computation formulas for certain classes of Apostol-type numbers and polynomials. The results of this paper contain not only the $\lambda$-Apostol-Daehee numbers and polynomials, but also Simsek numbers and polynomials, the Stirling numbers of the first kind, the Daehee numbers, and the Chu-Vandermonde identity. Furthermore, we derive an infinite series representation for the $\lambda$-Apostol-Daehee polynomials. By using functional equations containing the generating functions for the Cauchy numbers and the Riemann integrals of the generating functions for the $\lambda$-Apostol-Daehee numbers and polynomials, we also derive some identities and formulas for these numbers and polynomials. Moreover, we give implementation of a computation formula for the $\lambda$-Apostol-Daehee polynomials in Mathematica by Wolfram language. By this implementation, we also present some plots of these polynomials in order to investigate their behaviour some randomly selected special cases of its parameters. Finally, we conclude the paper with some comments and observations on our results.

\begin{quotation}
\bigskip
\end{quotation}

\noindent \textbf{2010 Mathematics Subject Classification.} 05A10, 05A15, 05A19, 11B37, 11B73, 11B83, 11S23, 33F05, 65D20.

\noindent \textbf{Key Words.} Generating functions, Functional Equations, Special numbers and polynomials, Stirling numbers of the first kind, Apostol-type numbers and polynomials, Simsek numbers and polynomials, Bernoulli numbers of the second kind, Daehee numbers and polynomials, Integral formulas, Chu-Vandermonde identity, Bernstein basis functions, Combinatorial sums, Mathematica implementation.

\section{Introduction}

In recent years, many studies on Apostol-type numbers and polynomials have been carried out by some researchers (see \cite{ChoiFAR}-\cite{KucukogluJNT2017}). Among others, in this paper, we are mainly dealt with the $\lambda$-Apostol-Daehee numbers $\mathfrak{D}_n\left(\lambda\right)$ and polynomials $\mathfrak{D}_n\left(x;\lambda\right)$ introduced and investigated by Simsek \cite{SimsekASCM2016, SimsekASCM2017} respectively as in the following generating functions:
\begin{equation}
	G_{\mathfrak{D}}\left(t; \lambda\right):=\frac{\log\lambda+\log \left(1+\lambda t\right)}{\lambda\left(1+\lambda t\right)-1}=\sum_{n=0}^{\infty}\mathfrak{D}_n\left(\lambda\right)\frac{t^n}{n!},
	\label{L-ADnum}
\end{equation}
and
\begin{equation}
	G_{\mathfrak{D}}\left(t, x; \lambda\right):=\frac{\log\lambda+\log \left(1+\lambda t\right)}{\lambda\left(1+\lambda t\right)-1}\left(1+\lambda t\right)^x=\sum_{n=0}^{\infty}\mathfrak{D}_n\left(x;\lambda\right)\frac{t^n}{n!},
	\label{L-ADpoly}
\end{equation}
(\textit{cf}. \cite{SimsekASCM2016, SimsekASCM2017}; and also see \cite{SimsekYardimci2016}).

First few values of the numbers $\mathcal{D}_{n}(\lambda )$ are given as follows:%
\begin{eqnarray*}
	\mathcal{D}_{0}(\lambda ) &=&\frac{\log \lambda }{\lambda -1},\\
	\mathcal{D}_{1}(\lambda ) &=&-\frac{\lambda ^{2}\log \lambda }{\left(
		\lambda -1\right) ^{2}}+\frac{\lambda }{\lambda -1}, \\
	\mathcal{D}_{2}(\lambda ) &=&\frac{2\lambda ^{4}\log \lambda }{\left(
		\lambda -1\right) ^{3}}+\frac{\lambda ^{2}\left( 1-3\lambda \right) }{\left(
		\lambda -1\right) ^{2}}, \\
	\mathcal{D}_{3}(\lambda ) &=&-\frac{6\lambda ^{6}\log \lambda }{\left(
		\lambda -1\right) ^{4}}+\frac{\lambda ^{3}\left( 11\lambda ^{2}-7\lambda
		+2\right) }{\left( \lambda -1\right) ^{3}},
\end{eqnarray*}
and so on (\textit{cf}. \cite{KucukogluAADM2019, SimsekASCM2016, SimsekASCM2017, SimsekYardimci2016}).

Another family of Apostol-type numbers and polynomials is the family of the numbers $Y_{n}\left( \lambda \right) $ (so-called Simsek numbers) and the polynomials $Y_{n}\left( x;\lambda \right) $ (so-called Simsek polynomials) defined respectively by the following generating functions (\textit{cf}. \cite{simsekTJM}):
\begin{equation}
	F\left( t,\lambda \right):=\frac{2}{\lambda \left( 1+\lambda t\right) -1}%
	=\sum\limits_{n=0}^{\infty }Y_{n}\left( \lambda \right) \frac{t^{n}}{n!},
	\label{YNumGenFunc}
\end{equation}%
and
\begin{equation}
	F\left( t,x,\lambda \right) :=\frac{2\left( 1+\lambda t\right) ^{x}}{\lambda
		\left( 1+\lambda t\right) -1}=\sum\limits_{n=0}^{\infty }Y_{n}\left(
	x;\lambda \right) \frac{t^{n}}{n!}.
	\label{YPolyGenFunc}
\end{equation}

For $n\in \mathbb{N}_{0}:=\{0,1,2,3,\dots\}$, the numbers $Y_{n}(\lambda )$ are computed by the following explicit formula:
\begin{equation}
	Y_{n}(\lambda )=2(-1)^{n}\frac{n!}{\lambda -1}\left( \frac{\lambda ^{2}}{%
		\lambda -1}\right) ^{n},
	\label{YnExpl}
\end{equation}
by which, one may easily compute first few values of the numbers $Y_{n}\left( \lambda \right) $ as below:
\begin{eqnarray*}
	&&Y_{0}(\lambda )=\frac{2}{\lambda -1}, \quad Y_{1}(\lambda )=-\frac{2\lambda ^{2}%
	}{\left( \lambda -1\right) ^{2}}, \quad Y_{2}(\lambda )=\frac{4\lambda ^{4}}{\left(
		\lambda -1\right) ^{3}}, \\
	&&Y_{3}(\lambda )=-\frac{12\lambda ^{6}}{\left( \lambda -1\right) ^{4}}%
	,\quad Y_{4}(\lambda )=\frac{48\lambda ^{8}}{\left( \lambda -1\right) ^{5}},
\end{eqnarray*}%
and so on (\textit{cf}. \cite{simsekTJM}; and also see \cite{KucukogluJNT2017}).

Let $\left( x\right) _{n}=x\left(x-1\right)\dots\left(x-n+1\right)$ with $\left( x\right) _{0}=1$. Then, the combination of (\ref{YNumGenFunc}) with (\ref{YPolyGenFunc}) yields the relation between the numbers $Y_{n}(\lambda )$ and the polynomials $Y_{n}(x;\lambda )$ given by 
\begin{equation}
	Y_{n}(x;\lambda )=\sum_{j=0}^{n}\binom{n}{j} Y_{j}(\lambda ) \lambda
	^{n-j}(x)_{n-j},
	\label{ttC}
\end{equation}
by which, one may easily compute first few values of the polynomials $Y_{n}\left(x; \lambda \right) $ as below:
\begin{eqnarray*}
	Y_{0}(x;\lambda ) &=&\frac{2}{\lambda -1}, \\
	Y_{1}(x;\lambda ) &=&\frac{2\lambda }{\lambda -1}x-\frac{2\lambda ^{2}}{%
		\left( \lambda -1\right) ^{2}}, \\
	Y_{2}(x;\lambda ) &=& \frac{2\lambda ^{2}}{\lambda -1}x^{2}-\frac{6\lambda
		^{3}-2\lambda ^{2}}{\left( \lambda -1\right) ^{2}}x+\frac{4\lambda ^{4}}{%
		\left( \lambda -1\right) ^{3}}, \\
	Y_{3}(x;\lambda ) &=&\frac{2\lambda ^{3}}{\lambda -1}x^{3}-\frac{12\lambda
		^{4}-6\lambda ^{3}}{\left( \lambda -1\right) ^{2}}x^{2}+\frac{22\lambda
		^{5}-14\lambda ^{4}+4\lambda ^{3}}{\left( \lambda -1\right) ^{3}}x-\frac{%
		12\lambda ^{6}}{\left( \lambda -1\right) ^{4}},
\end{eqnarray*}
and so on (\textit{cf}. \cite{simsekTJM}; and also see \cite{KucukogluJNT2017}).

Another family of Apostol-type numbers and polynomials is the family of the numbers $Y_{n}^{\left( -k\right) }\left(
\lambda \right)$ (so-called negative higher-order Simsek numbers) and the polynomials $Q_{n}\left(x; \lambda,
k\right)$ (so-called negative higher-order Simsek polynomials) defined respectively by the following generating functions:  
\begin{equation}
	\mathcal{G}_{Y}\left(t,k;\lambda\right):=2^{-k}\left(\lambda \left( 1+\lambda
	t\right) -1\right) ^{k}=\sum_{n=0}^{\infty }Y_{n}^{\left( -k\right) }\left(
	\lambda \right) \frac{t^{n}}{n!}  \label{GenFHigOrdNegYpoly}
\end{equation}%
and
\begin{equation}
	\mathcal{G}_{Q}\left(t,x,k;\lambda\right):=\mathcal{G}_{Y}\left(t,k;\lambda\right)
	\left(1+\lambda t\right)^x=\sum_{n=0}^{\infty }Q_{n}\left(x; \lambda,
	k\right) \frac{t^{n}}{n!}, \label{GenFHigOrdNegYpolyx}
\end{equation}
(\textit{cf}. \cite{KucukogluAxiom2019}).

For $n\in \mathbb{N}_{0}$, the numbers $Y_{n}^{\left( -k\right) }\left(\lambda \right)$ are computed by the following explicit formula:
\begin{equation}
	Y_{n}^{\left( -k\right) }\left( \lambda \right)=\left\{ 
	\begin{array}{cc}
		2^{-k} n!\binom{k}{n}\lambda^{2n}\left(\lambda-1\right)^{k-n} & \text{if}
		\quad n\leq k \\ 
		0 & \text{if} \quad n>k
	\end{array}
	\right.  \label{HigOrdNegYpolyExplicit}
\end{equation}
by which, one may easily compute the values of the numbers $Y_{n}^{\left( -k\right) }\left(\lambda \right)$ as below:
\begin{eqnarray*}
	Y_{0}^{\left( -k\right) }\left( \lambda
	\right)&=&2^{-k}\left(\lambda-1\right)^{k}, \\
	Y_{1}^{\left( -k\right) }\left( \lambda \right)&=&2^{-k} \binom{k}{1}%
	\lambda^{2}\left(\lambda-1\right)^{k-1}, \\
	Y_{2}^{\left( -k\right) }\left( \lambda \right)&=&2^{-k} 2!\binom{k}{2}%
	\lambda^{4}\left(\lambda-1\right)^{k-2}, \\
	&\vdots& \\
	Y_{j}^{\left( -k\right) }\left( \lambda \right)&=&2^{-k} j!\binom{k}{j}%
	\lambda^{2j}\left(\lambda-1\right)^{k-j} \quad \mathit{for} \quad j\leq k, \\
	&\vdots& \\
	Y_{k}^{\left( -k\right) }\left( \lambda \right)&=&2^{-k} k!\lambda^{2k},
\end{eqnarray*}
(\textit{cf}. \cite{KucukogluAxiom2019}).

The combination of (\ref{GenFHigOrdNegYpoly}) with (\ref{GenFHigOrdNegYpolyx}) yields the relation between the numbers $Y_{n}^{\left( -k\right) }\left(\lambda \right)$ and the polynomials $Q_{n}\left(x; \lambda, k\right)$ is given, for $k, n \in \mathbb{N}_0$, by
\begin{equation}
	Q_{n}\left(x; \lambda, k\right)=\sum_{j=0}^{n}\binom{n}{j}Y_{j}^{\left( -k\right) }\left( \lambda
	\right) %
	\lambda^{n-j}\left(x\right)_{n-j},  \label{Formula-HigOrdNegYpolyx}
\end{equation}
by which, one may easily compute first few values of the polynomials $Q_{n}\left(x; \lambda, k\right)$ as below:
\begin{eqnarray*}
	Q_{0}\left(x; \lambda, k\right)&=&2^{-k}\left(\lambda-1\right)^k,  \notag \\
	Q_{1}\left(x; \lambda, k\right)&=&2^{-k}\left(\lambda-1\right)^k\lambda
	x+2^{-k}k\lambda^2\left(\lambda-1\right)^{k-1},  \notag \\
	Q_{2}\left(x; \lambda,
	k\right)&=&2^{-k}\left(\lambda-1\right)^k\lambda^2x^2+\left(-2^{-k}\left(%
	\lambda-1\right)^k\lambda^2
	+2^{-k+1}k\lambda^3\left(\lambda-1\right)^{k-1}\right)x  \notag \\
	&&+2^{-k}k\left(k-1\right)\lambda^4\left(\lambda-1\right)^{k-1},
\end{eqnarray*}
and so on (\textit{cf}. \cite{KucukogluAxiom2019}).

Let $x\in \left[ 0,1\right] $ and $k \in \mathbb{N}_0$. Then, the Bernstein Basis functions, $B_{k}^{n}(x)$, is defined as below:
\begin{equation}
	B_{k}^{n}(x)=\binom{n}{k} x^{k}(1-x)^{n-k}, \qquad \left(k=0,1,\dots,n; \ \ n \in \mathbb{N}_0\right)
	\label{Berns-Basis}
\end{equation}
and its generating function is given by 
\begin{equation}
	\frac{\left(xt\right)^{k}e^{\left( 1-x\right) t}}{k!}%
	=\sum_{n=0}^{\infty}B_{k}^{n}(x)\frac{t^{n}}{n!},  \label{GenFunc-Berns}
\end{equation}
so that the Bernstein Basis functions have relationships with a large number of concepts including the Bezier curves, the binomial distribution, the Poisson distribution, the Catalan numbers, and etc.; see, for details, \cite{AcikgozSerkan, ErkusDuman, Duman, Lorentz,FPTASimsek2013,SimsekBVP2013, SimsekHJMS2014, SimsekMMAS2015, Simsek Acikgoz} and also cited references therein.

It is concluded with the help of (\ref{HigOrdNegYpolyExplicit}) and (\ref{Berns-Basis}) that there exists the following relation between the numbers $Y_{n}^{\left( -k\right) }\left( \lambda \right) $ and the Bernstein basis functions:  
\begin{equation}
	Y_{n}^{\left( -k\right) }\left( \lambda \right)=\frac{\left(-1\right)^{k-n} n!}{2^{k}}\lambda^n B_{n}^{k}(\lambda)
	\label{Y-Berns}
\end{equation}
where $n, k \in \mathbb{N}_0$ and $\lambda\in \left[ 0,1\right] $ (\textit{cf}. \cite{KucukogluAxiom2019}).

Actually, the numbers $Y_{n}^{\left( -k\right) }\left( \lambda \right)$ have other relations than its relation to the Bernstein basis functions. Among others, the numbers $Y_{n}^{\left( -k\right) }\left( \lambda \right)$ have relationships with the Poisson--Charlier polynomials, the Bell polynomials (i.e., exponential polynomials) and other kinds of combinatorial numbers. To see the relations mentioned above, the interested reader may glance at the recent paper \textup{\cite{KucukogluAxiom2019}}.

The Stirling numbers of the first kind, $S_{1}(n,k)$, are defined, for $n, k\in \mathbb{N}_0$, by the following recurrence relation:
\begin{equation*}
	S_{1}(n+1,k)=-nS_{1}(n,k)+S_{1}(n,k-1)
\end{equation*}%
with the side conditions $S_{1}(0,0)=1$, $S_{1}(0,k)=0$ if $k>0$, $S_{1}(n,0)=0$ if $n>0$, $%
S_{1}(n,k)=0$ if $k>n$; and these numbers are also given by
\begin{equation}
	\left( x\right) _{n}=\sum_{k=0}^{n}S_{1}\left( n,k\right) x^{k}
	\label{DefinitionFirstStirling}
\end{equation}
and 
\begin{equation}
	\frac{\left( \log (1+t)\right) ^{k}}{k!}=\sum_{n=k}^{\infty }S_{1}(n,k)\frac{%
		t^{n}}{n!}
	\label{S1}
\end{equation}
(\textit{cf}. \cite{BonaBOOK, CharalambidesBOOK, Comtet, Gradimir, Qi, SimsekFPTA}; and the references cited therein).

The Cauchy numbers (or the Bernoulli numbers of the second kind), $b_{n} \left(0\right)$, are defined by
\begin{equation}
	b_{n} \left(0\right) =\int\limits_{0}^{1}\left( x\right) _{n}dx 
	\label{Cauchy-Int-Formula}
\end{equation}
and
\begin{equation}
	G_{C}\left(t\right):=\frac{t}{\log \left(t+1\right) }=\sum_{n=0}^{%
		\infty }b_{n} \left(0\right)\frac{t^{n}}{n!}  \label{Cauchy-1}
\end{equation}%
(\textit{cf}. \cite[p. 116]{Roman}, \cite{KimEtAl2016}, \cite{MerliniCauchy}, \cite{Qi}). 

The combination of (\ref{DefinitionFirstStirling}) with (\ref{Cauchy-Int-Formula}) yields the relation between the Cauchy numbers and the Stirling numbers of the first kind given by
\begin{equation}
	b_{n} \left(0\right)= \sum\limits_{m=0}^{n}\frac{%
		S_{1}\left(n,m\right) }{m+1},
	\label{Cauchy-Int-Formula2}
\end{equation}
(\textit{cf}. \cite[p. 294]{Comtet}, \cite[p. 1908]{MerliniCauchy}, \cite[p. 114]{Roman}).

The Daehee numbers, $D_n$, are defined by the following generating function:
\begin{equation}
	F_D\left(t\right):=\frac{\log\left(1+t\right)}{t}=\sum_{n=0}^{\infty}D_n \frac{t^n}{n!}
	\label{DaeheeNumGF}
\end{equation}
and these numbers are computed, for $n\in \mathbb{N}_{0}$, by the following explicit formula:
\begin{equation}
	D_n=\frac{\left(-1\right)^n n!}{n+1}
	\label{DaeheeNumExpl}
\end{equation}
(\textit{cf}. \cite{El-Desouky, KimDaehee, SimsekASCM2016, SimsekASCM2017}).

The well-known Chu-Vandermonde identity is given by
\begin{equation}
	\binom{x+y}{n}=\frac{1}{n!}
	\sum\limits_{k=0}^{n}\binom{n}{k}\left(x\right)_{k}\left(y\right)_{n-k}
	\label{Chu-Vandermonde}
\end{equation}
(\textit{cf}. \cite{Comtet, Jordan, ysimsekMTJPAM}).

The outline of this paper may briefly given as follows: 

In Section \ref{Section-2}, we present various identities and computation formulas containing not only the $\lambda$-Apostol-Daehee numbers and polynomials, but also Simsek numbers and polynomials, the Stirling numbers of the first kind, the Daehee numbers and the Chu-Vandermonde identity. Besides, we derive an infinite series representation for the $\lambda$-Apostol-Daehee polynomials.

In Section \ref{Section-3}, by using functional equations containing the generating functions for the Cauchy numbers and the integrals of the generating functions for the $\lambda$-Apostol-Daehee numbers and polynomials, we also derive some identities and formulas for these numbers and polynomials.

In Section \ref{Section-4}, we give Mathematica implementation of a formula which computes the $\lambda$-Apostol-Daehee polynomials in terms of the Simsek polynomials. By this implementation, some plots of the $\lambda$-Apostol-Daehee polynomials are presented for some randomly selected special cases.

In Section \ref{Section-5}, we conclude the paper with some comments and observations on our results.

\section{Identities containing Apostol-type numbers and polynomials}
\label{Section-2}

In this section, by using the techniques of generating function and their functional equations, we derive some identities involving not only the Chu-Vandermonde identity, but also some special numbers and polynomials such as the numbers $\mathfrak{D}_n\left(\lambda\right)$, the polynomials $\mathfrak{D}_n\left(x;\lambda\right)$, the numbers $Y_n\left(\lambda\right)$, the polynomials $Y_n\left(x;\lambda\right)$, the numbers $S_1\left(n,m\right)$ and the numbers $D_n$. In addition, we get computation formulas for not only the numbers $\mathfrak{D}_n\left(\lambda\right)$, but also the polynomials $\mathfrak{D}_n\left(x;\lambda\right)$. Moreover, we derive an infinite series representation for the polynomials $\mathfrak{D}_m\left(x;\lambda\right)$ in terms of the polynomials $Q_{m}\left(x; \lambda, n\right)$.

By (\ref{L-ADpoly}), we have
\begin{equation*}
	\left(1+\lambda t\right)^x=\frac{\lambda\left(1+\lambda t\right)-1}{\log\lambda+\log \left(1+\lambda t\right)}\sum_{n=0}^{\infty}\mathfrak{D}_n\left(x;\lambda\right)\frac{t^n}{n!}
\end{equation*}
and
\begin{equation*}
	\left(1+\lambda t\right)^y=\frac{\lambda\left(1+\lambda t\right)-1}{\log\lambda+\log \left(1+\lambda t\right)}\sum_{n=0}^{\infty}\mathfrak{D}_n\left(y;\lambda\right)\frac{t^n}{n!}.
\end{equation*}
Multiplying the above two equations each other, we get
\begin{equation*}
	\left(1+\lambda t\right)^{x+y}=\left(\frac{\lambda\left(1+\lambda t\right)-1}{\log\lambda+\log \left(1+\lambda t\right)}\right)^2 \sum_{n=0}^{\infty}\mathfrak{D}_n\left(x;\lambda\right)\frac{t^n}{n!}\sum_{n=0}^{\infty}\mathfrak{D}_n\left(y;\lambda\right)\frac{t^n}{n!}.
\end{equation*}
With the application of the Binomial theorem and the Cauchy product rule to the above equation, we get
\begin{equation*}
	\sum_{n=0}^{\infty}\binom{x+y}{n}\lambda^n t^n=\frac{\left(\lambda\left(1+\lambda t\right)-1\right)^2}{\left(\log\lambda\right)^2\left(1+\frac{\log \left(1+\lambda t\right)}{\log\lambda}\right)^2} \sum_{n=0}^{\infty}\sum_{j=0}^{n}\binom{n}{j}\mathfrak{D}_j\left(x;\lambda\right)\mathfrak{D}_{n-j}\left(y;\lambda\right)\frac{t^n}{n!}.
\end{equation*}
Thus, we have
\begin{eqnarray*}
	\sum_{n=0}^{\infty}\binom{x+y}{n}\lambda^n t^n&=&
	\frac{\left(\lambda\left(1+\lambda t\right)-1\right)^2}{\left(\log\lambda\right)^2}
	\sum_{n=0}^{\infty}\binom{-2}{n}\frac{\left(\log\left(1+\lambda t\right)\right)^n}{\left(\log \lambda\right)^n}\\
	&&\times\sum_{n=0}^{\infty}\sum_{j=0}^{n}\binom{n}{j}\mathfrak{D}_j\left(x;\lambda\right)\mathfrak{D}_{n-j}\left(y;\lambda\right)\frac{t^n}{n!}.
\end{eqnarray*}
By combining (\ref{S1}) with the above equation, after some elementary calculations, we get
\begin{eqnarray*}
	\sum_{n=0}^{\infty}\binom{x+y}{n}\lambda^n t^n&=&
	\frac{\left(\lambda\left(1+\lambda t\right)-1\right)^2}{\left(\log\lambda\right)^2}
	\sum_{m=0}^{\infty}\sum_{n=0}^{m}\binom{-2}{n}n!\frac{S_1\left(m,n\right)}{\left(\log \lambda\right)^n}\frac{t^m}{m!}\\
	&&\times\sum_{n=0}^{\infty}\sum_{j=0}^{n}\binom{n}{j}\mathfrak{D}_j\left(x;\lambda\right)\mathfrak{D}_{n-j}\left(y;\lambda\right)\frac{t^n}{n!}.
\end{eqnarray*}
By applying the Cauchy product rule to the above equation, we get
\begin{eqnarray*}
	\sum_{m=0}^{\infty}\binom{x+y}{m}\lambda^m t^m&=&
	\frac{\left(\lambda^4 t^2+2\lambda^2\left(\lambda-1\right)t+\left(\lambda-1\right)^2\right)}{\left(\log\lambda\right)^2}\\
	&&\times
	\sum_{m=0}^{\infty}\sum_{k=0}^{m}\binom{m}{k}\sum_{n=0}^{k}\binom{-2}{n}n!\frac{S_1\left(k,n\right)}{\left(\log \lambda\right)^n}\\
	&&\times\sum_{j=0}^{m-k}\binom{m-k}{j}\mathfrak{D}_j\left(x;\lambda\right)\mathfrak{D}_{m-k-j}\left(y;\lambda\right)\frac{t^m}{m!}.
\end{eqnarray*}
After some elementary calculations and by comparing the coefficients of $\frac{t^m}{m!}$ on both sides of the above equation, we arrive at the following theorem:
\begin{theorem}
	Let $m \in \mathbb{N}\setminus\{0,1\}$. Then we have
	\begin{eqnarray}
		\binom{x+y}{m}&=&\frac{\lambda^{4-m}}{\left(\log\lambda\right)^2}m\left(m-1\right)A\left(m-2\right) \label{Th-1a}\\
		&&+\frac{2\lambda^{2-m}\left(\lambda-1\right)}{\left(\log\lambda\right)^2}mA\left(m-1\right) \notag\\
		&&+\frac{\lambda^{-m}\left(\lambda-1\right)^2}{\left(\log\lambda\right)^2}A\left(m\right),\notag
	\end{eqnarray}
	where
	\begin{eqnarray*}
		A\left(m\right)=\sum_{k=0}^{m}\sum_{n=0}^{k}\sum_{j=0}^{m-k}\binom{m}{k}\binom{-2}{n}\binom{m-k}{j}\frac{n!S_1\left(k,n\right)\mathfrak{D}_j\left(x;\lambda\right)\mathfrak{D}_{m-k-j}\left(y;\lambda\right)}{\left(\log \lambda\right)^n}.
	\end{eqnarray*}
\end{theorem}

Combining (\ref{Th-1a}) with (\ref{Chu-Vandermonde}) yields the following corollary:
\begin{corollary}
	Let $m \in \mathbb{N}\setminus\{0,1\}$. Then we have
	\begin{eqnarray}
		\sum\limits_{k=0}^{m}\binom{m}{k}\left(x\right)_{k}\left(y\right)_{m-k} &=&\frac{\lambda^{4-m}}{\left(\log\lambda\right)^2}m\left(m-1\right)A\left(m-2\right) \label{Th-1}\\
		&&+\frac{2\lambda^{2-m}\left(\lambda-1\right)}{\left(\log\lambda\right)^2}mA\left(m-1\right) \notag\\
		&&+\frac{\lambda^{-m}\left(\lambda-1\right)^2}{\left(\log\lambda\right)^2}A\left(m\right).\notag
	\end{eqnarray}
\end{corollary}

By the combination of (\ref{L-ADpoly}) with (\ref{YPolyGenFunc}) and (\ref{DaeheeNumGF}), we get the following functional equation:
\begin{equation}
	G_{\mathfrak{D}}\left(t, x; \lambda\right)=\left(\frac{\log\lambda}{2}+\frac{\lambda t F_D\left(\lambda t\right)}{2}\right)F\left( t,x,\lambda \right).
\end{equation}
which yields
\begin{equation*}
	\sum_{n=0}^{\infty}\mathfrak{D}_n\left(x;\lambda\right)\frac{t^n}{n!}=\left(\frac{\log\lambda}{2}+\frac{\lambda t}{2}\sum_{n=0}^{\infty}\lambda^n D_n \frac{t^n}{n!}\right)\sum_{n=0}^{\infty} {Y}_n\left(x;\lambda\right)\frac{t^n}{n!}.
\end{equation*}
By applying the Cauchy product rule to the above equation, after some elementary calculations, we get
\begin{eqnarray*}
	\sum_{n=0}^{\infty}\mathfrak{D}_n\left(x;\lambda\right)\frac{t^n}{n!}&=&\frac{\log\lambda}{2}\sum_{n=0}^{\infty} {Y}_n\left(x;\lambda\right)\frac{t^n}{n!}\\
	&&+\frac{1}{2}\sum_{n=0}^{\infty}\sum_{j=0}^{n-1}n\binom{n-1}{j}\lambda^{n-j} D_{n-j-1}{Y}_j\left(x;\lambda\right)\frac{t^n}{n!}.
\end{eqnarray*}
Comparing the coefficients of $\frac{t^n}{n!}$ on both sides of the above equation yields the following theorem:
\begin{theorem}
	Let $n \in \mathbb{N}$. Then we have
	\begin{equation}
		\mathfrak{D}_n\left(x;\lambda\right)=\frac{\log\lambda}{2}{Y}_n\left(x;\lambda\right)+\frac{1}{2}\sum\limits_{j=0}^{n-1}n\binom{n-1}{j}\lambda^{n-j} D_{n-j-1}{Y}_j\left(x;\lambda\right).
		\label{Th-2}
	\end{equation}
\end{theorem}

Using (\ref{DaeheeNumExpl}) in (\ref{Th-2}), we get a relation between $\lambda$-Apostol-Daehee polynomials and Simsek polynomials, given by following corollary:
\begin{corollary}
	Let $n \in \mathbb{N}$. Then we have
	\begin{equation}
		\mathfrak{D}_n\left(x;\lambda\right)=\frac{\log\lambda}{2}{Y}_n\left(x;\lambda\right)-\frac{n!}{2}\sum\limits_{j=0}^{n-1}\left(-1\right)^{n-j}\frac{\lambda^{n-j} {Y}_j\left(x;\lambda\right)}{j!\left(n-j\right)}.
		\label{Cor1}
	\end{equation}
\end{corollary}

Substituting $x=0$ into (\ref{Cor1}), we also get a relation, between $\lambda$-Apostol-Daehee numbers and Simsek numbers, given by following corollary:
\begin{corollary}
	Let $n \in \mathbb{N}$. Then we have
	\begin{equation}
		\mathfrak{D}_n\left(\lambda\right)=\frac{\log\lambda}{2}{Y}_n\left(\lambda\right)+\frac{n!}{2}\sum\limits_{j=0}^{n-1}\left(-1\right)^{n-j-1}\frac{\lambda^{n-j}{Y}_j\left(\lambda\right)}{j!\left(n-j\right)}.
		\label{Cor2}
	\end{equation}
\end{corollary}

Combining (\ref{YnExpl}) with (\ref{Cor2}), we get a computation formula for the numbers $\mathfrak{D}_n\left(\lambda\right)$ by the following corollary:
\begin{corollary}
	\begin{equation}
		\mathfrak{D}_n\left(\lambda\right)=\frac{(-1)^{n}n!}{\lambda-1}\left( \left( \frac{\lambda ^{2}}{%
			\lambda -1}\right) ^{n} \log\lambda-\lambda^n\sum\limits_{j=0}^{n-1}\frac{1}{n-j}\left(\frac{\lambda}{\lambda-1}\right)^j\right).
		\label{Cor3}
	\end{equation}
\end{corollary}

\begin{remark}
	The computation formula \textup{(\ref{Cor3})}, obtained by reduction from the equation \textup{(\ref{Th-2})}, may also be obtained with the help of the application of the binomial theorem on the generating function for the numbers $\mathfrak{D}_n\left(\lambda\right)$. In the meanwhile, for another form of this formula, the interested reader may see the recent paper \textup{\cite[Theorem 8, p. 492]{KucukogluAADM2019}} in which other methods and generating function families used in order to achieve the aforementioned formula.
\end{remark}

By (\ref{Cor3}), we obtain a finite sum whose values is computed by the numbers $\mathfrak{D}_n\left(\lambda\right)$ given by the following corollary:
\begin{corollary}
	Let $n \in \mathbb{N}$. Then we have
	\begin{equation}
		\sum\limits_{j=0}^{n-1}\frac{1}{n-j}\left(\frac{\lambda}{\lambda-1}\right)^j=(-1)^{n+1}\frac{\left(\lambda-1\right)\mathfrak{D}_n\left(\lambda\right)}{n!\lambda^n}+\left( \frac{\lambda}{%
			\lambda -1}\right)^{n} \log\lambda.
		\label{Cor4}
	\end{equation}
\end{corollary}

By using the Taylor series expansion of the function $\log\left(1+\left(\lambda\left(1+\lambda t\right)-1\right)\right)$, by assuming that $| \lambda\left(1+\lambda t\right)-1 |<1$, in the equation (\ref{L-ADpoly}), and then by making some simplifications, we get
\begin{equation*}
	\sum_{m=0}^{\infty}\mathfrak{D}_m\left(x;\lambda\right)\frac{t^m}{m!}=\left(1+\lambda t\right)^x\sum_{n=0}^{\infty}\left(-1\right)^n\frac{\left(\lambda\left(1+\lambda t\right)-1\right)^{n}}{n+1}.
\end{equation*}
By combining (\ref{GenFHigOrdNegYpolyx}) with the above equation, we get
\begin{equation*}
	\sum_{m=0}^{\infty}\mathfrak{D}_m\left(x;\lambda\right)\frac{t^m}{m!}=\sum_{n=0}^{\infty}\left(-1\right)^n \frac{2^n}{n+1}\sum_{m=0}^{\infty }Q_{m}\left(x; \lambda,
	n\right) \frac{t^{m}}{m!}.
\end{equation*}
which yields
\begin{equation*}
	\sum_{m=0}^{\infty}\mathfrak{D}_m\left(x;\lambda\right)\frac{t^m}{m!}=\sum_{n=0}^{\infty }\sum_{m=0}^{\infty}\left(-1\right)^n \frac{2^n Q_{m}\left(x; \lambda,
		n\right)}{n+1} \frac{t^{m}}{m!}.
\end{equation*}
By assuming that $| \lambda-1 |<1$ and comparing the coefficients of $\frac{t^{m}}{m!}$ on both sides of the above equation yields a relation, between the numbers $\mathfrak{D}_m\left(x;\lambda\right)$ and the polynomials $Q_{m}\left(x; \lambda,
n\right)$, given by the following theorem:
\begin{theorem}
	If $| \lambda-1 |<1$, then we have the following infinite series representation for the polynomials $\mathfrak{D}_m\left(x;\lambda\right)$:
	\begin{equation}
		\mathfrak{D}_m\left(x;\lambda\right)=\sum_{n=0}^{\infty}\left(-1\right)^n \frac{2^n Q_{m}\left(x; \lambda,
			n\right)}{n+1}.
		\label{Th-Series-DQ}
	\end{equation}
\end{theorem}

\section{Further identities derived from integral formulas and Cauchy numbers}
\label{Section-3}

In this section, by using functional equations involving the generating functions for the Cauchy numbers and the integrals of the generating functions for the numbers $\mathfrak{D}_n\left(\lambda\right)$, the polynomials $\mathfrak{D}_n\left(x; \lambda\right)$, we derive some identities and formulas.

Integrating both-sides of the equation (\ref{L-ADpoly}), with respect to the variable $x$, from $0$ to $1$, we get the following integral formula:
\begin{eqnarray}
	\int\limits_{0}^{1} G_{\mathfrak{D}}\left(t, x; \lambda\right)\mathrm{d}x&=&\int\limits_{0}^{1} \frac{\log\lambda+\log \left(1+\lambda t\right)}{\lambda\left(1+\lambda t\right)-1}\left(1+\lambda t\right)^x \mathrm{d}x\\
	&=&\frac{\lambda t\left(\log\lambda+\log \left(1+\lambda t\right)\right)}{\left(\lambda\left(1+\lambda t\right)-1\right)\log\left(1+\lambda t\right)},
\end{eqnarray}
which, by (\ref{Cauchy-1}) and (\ref{L-ADnum}), yields the following functional equation:
\begin{eqnarray}
	\int\limits_{0}^{1} G_{\mathfrak{D}}\left(t, x; \lambda\right)\mathrm{d}x=G_{\mathfrak{D}}\left(t; \lambda\right)G_{C}\left(\lambda t\right).
\end{eqnarray}
Combining the above equation with (\ref{Cauchy-1}), (\ref{L-ADnum}) and (\ref{L-ADpoly}) yields
\begin{eqnarray}
	\int\limits_{0}^{1} \sum_{n=0}^{\infty}\mathfrak{D}_n\left(x;\lambda\right)\frac{t^n}{n!}\mathrm{d}x=\sum_{n=0}^{\infty}\mathfrak{D}_n\left(\lambda\right)\frac{t^n}{n!}\sum_{n=0}^{%
		\infty }\lambda^n b_{n} \left(0\right)\frac{t^{n}}{n!}.
\end{eqnarray}
By applying the Cauchy product rule to the right-hand side of the above equation, after some elementary calculations, we get
\begin{eqnarray}
	\sum_{n=0}^{\infty}\int\limits_{0}^{1}\mathfrak{D}_n\left(x;\lambda\right)\mathrm{d}x\frac{t^n}{n!}=\sum_{n=0}^{\infty}\sum_{m=0}^{n}\binom{n}{m}\lambda^m b_{m} \left(0\right)\mathfrak{D}_{n-m}\left(\lambda\right)\frac{t^{n}}{n!}.
\end{eqnarray}
Comparing the coefficients of $\frac{t^{n}}{n!}$ on both sides of the above equation yields the following theorem:
\begin{theorem}
	\begin{eqnarray}
		\int\limits_{0}^{1}\mathfrak{D}_n\left(x;\lambda\right)\mathrm{d}x=\sum_{m=0}^{n}\binom{n}{m}\lambda^m b_{m} \left(0\right)\mathfrak{D}_{n-m}\left(\lambda\right).
		\label{Th-Cauchy}
	\end{eqnarray}
\end{theorem}

\begin{remark}
	By using the generating function for the $k$-th order $\lambda$-Apostol-Daehee polynomials, Choi \textup{\cite[Theorem 5, p.1854]{ChoiFAR}} gave the following integral formula:
	\begin{equation*}
		\int\limits_{\alpha}^{\alpha+1}\mathfrak{D}^{\left(k\right)}_n\left(x;\lambda\right)\mathrm{d}x=\sum_{m=0}^{n}m!\binom{n}{m}\lambda^{m}\mathfrak{D}^{\left(k\right)}_{n-m}\left(\alpha;\lambda\right)p_m.
	\end{equation*}
	If we substitute $k=1$ and $\alpha=0$ into the above formula, we get
	\begin{equation*}
		\int\limits_{0}^{1}\mathfrak{D}_n\left(x;\lambda\right)\mathrm{d}x=\sum_{m=0}^{n}m!\binom{n}{m}\lambda^{m}\mathfrak{D}_{n-m}\left(\lambda\right)p_m.
	\end{equation*}
	When we compare the above formula with the equation \textup{(\ref{Th-Cauchy})}, it is shown that the numbers $m!p_m$ considered in the formula above actually correspond to the Cauchy numbers $b_{m} \left(0\right)$ which is obtained by the techniques of generating functions and their functional equations. Thus, we conclude that Choi \textup{\cite{ChoiFAR}} modified the numbers $b_{m} \left(0\right)$ as follows:
	\begin{equation*}
		m!p_m=b_{m} \left(0\right)
	\end{equation*}
	in order to obtain an integral formula for the higher-order $\lambda$-Apostol-Daehee polynomials.
\end{remark}

Integrating both-sides of the equation (\ref{L-ADpoly}), with respect to the variable $x$, from $0$ to $z$, we get the following integral formula:
\begin{eqnarray}
	\int\limits_{0}^{z} G_{\mathfrak{D}}\left(t, x; \lambda\right)\mathrm{d}x&=&\int\limits_{0}^{z} \frac{\log\lambda+\log \left(1+\lambda t\right)}{\lambda\left(1+\lambda t\right)-1}\left(1+\lambda t\right)^x \mathrm{d}x\\
	&=&\frac{\left(\left(1+\lambda t\right)^z-1\right)\left(\log\lambda+\log \left(1+\lambda t\right)\right)}{\left(\lambda\left(1+\lambda t\right)-1\right)\log\left(1+\lambda t\right)}
\end{eqnarray}
which, by (\ref{L-ADnum}) and (\ref{L-ADpoly}), yields the following functional equation:
\begin{eqnarray}
	\int\limits_{0}^{z} G_{\mathfrak{D}}\left(t, x; \lambda\right)\mathrm{d}x=\frac{\left( G_{\mathfrak{D}}\left(t, z; \lambda\right)- G_{\mathfrak{D}}\left(t; \lambda\right)\right)G_{C}\left(\lambda t\right)}{\lambda t}.
\end{eqnarray}
Combining the above equation with (\ref{Cauchy-1}), (\ref{L-ADnum}) and (\ref{L-ADpoly}) yields
\begin{eqnarray*}
	\int\limits_{0}^{z} \sum_{n=0}^{\infty}\mathfrak{D}_n\left(x;\lambda\right)\frac{t^n}{n!}\mathrm{d}x=\frac{1}{\lambda t}\left( \sum_{n=0}^{\infty}\left(\mathfrak{D}_n\left(z;\lambda\right)-\mathfrak{D}_n\left(\lambda\right)\right)\frac{t^n}{n!}\right)\sum_{n=0}^{\infty }\lambda^n b_{n} \left(0\right)\frac{t^{n}}{n!}.
\end{eqnarray*}
By applying the Cauchy product rule to the right-hand side of the above equation, after some elementary calculations, we get
\begin{eqnarray*}
	\int\limits_{0}^{z} \sum_{n=0}^{\infty}\mathfrak{D}_n\left(x;\lambda\right)\frac{t^n}{n!}\mathrm{d}x&=&\sum_{n=0}^{\infty}\frac{1}{n+1}\sum_{m=0}^{n+1}\binom{n+1}{m}\lambda^{m-1} b_{m} \left(0\right)\\
	&&\times\left(\mathfrak{D}_{n+1-m}\left(z;\lambda\right)-\mathfrak{D}_{n+1-m}\left(\lambda\right)\right)\frac{t^{n}}{n!}.
\end{eqnarray*}
Comparing the coefficients of $\frac{t^{n}}{n!}$ on both sides of the above equation yields the following theorem:
\begin{theorem}
	\begin{eqnarray}
		&&\int\limits_{0}^{z}\mathfrak{D}_n\left(x;\lambda\right)\mathrm{d}x\\
		&&=\frac{1}{n+1}\sum_{m=0}^{n+1}\binom{n+1}{m}\lambda^{m-1} b_{m} \left(0\right)\left(\mathfrak{D}_{n+1-m}\left(z;\lambda\right)-\mathfrak{D}_{n+1-m}\left(\lambda\right)\right).\notag
	\end{eqnarray}
\end{theorem}

\section{Implementation of computation formulas involving $\lambda$-Apostol-Daehee polynomials}
\label{Section-4}

In this section, by implementing some of our results with the aid of the Wolfram programming language in Mathematica \cite{WolframCloud},  we compute a few values of the $\lambda$-Apostol-Daehee polynomials. In addition, we also give some illustrations involving two dimensional plots of the $\lambda$-Apostol-Daehee polynomials.

We first give Mathematica implementation of the equation (\ref{Th-2}) in Implementation \ref{YlistDpoly} in which we utilized from the Implementation \ref{YlistYpoly} and the Implementation \ref{YlistNum} given by Simsek and Kucukoglu \cite{SimsekKucukoglu2020Chapter} in order to compute the rational functions $Y_{n}\left(\lambda \right)$ and the polynomials $Y_{n}\left(x; \lambda \right)$.

\begin{lstlisting}[language=Mathematica, label=YlistDpoly, caption={Let the expression \texttt{lparameter} denote the parameter $\lambda$ of the polynomials $\mathfrak{D}_n\left(x;\lambda \right)$. Then, the following Mathematica code returns the values of these polynomials}.]
DPoly[x_,lparameter_,n_]:=(Log[lparameter]/2)*YPoly[x,lparameter,n]+(Factorial[n]/2)*Sum[((-1)^(n-j-1))*((((lparameter)^(n-j))*YPoly[x,lparameter,j])/(Factorial[j]*(n-j))), {j,0,n-1}]
\end{lstlisting}

\begin{lstlisting}[language=Mathematica, label=YlistYpoly, caption={Let the expression \texttt{lparameter} denote the parameter $\lambda$ of the polynomials $Y_{n}\left(x; \lambda \right)$. Then, the following Mathematica code returns the values of  these polynomials} (\textit{cf}. \cite{SimsekKucukoglu2020Chapter}).]
YPoly[x_,lparameter_,n_]:=Sum[Binomial[n,j]*(lparameter^(n-j))*FactorialPower[x, n-j, 1]*YNum[lparameter,j], {j,0,n}]
\end{lstlisting}

\begin{lstlisting}[language=Mathematica, label=YlistNum, caption={Let the expression \texttt{lparameter} denote the parameter $\lambda$. Then, the following Mathematica code returns the values of the rational functions $Y_{n}\left(\lambda \right)$} (\textit{cf}. \cite{SimsekKucukoglu2020Chapter}).]
YNum[l_,n_]:=2*((-1)^n)* (Factorial[n]/(l-1))*(((l^2)/(l-1))^n)
\end{lstlisting}

By using the Implementation \ref{YlistDpoly} in Mathematica, we compute first three values of the $\lambda$-Apostol-Daehee polynomials as follows:
\begin{eqnarray*}
	&&\mathfrak{D}_0\left(x;\lambda\right)=\frac{\log \lambda}{\lambda -1}, \\
	&&\mathfrak{D}_1\left(x;\lambda\right)=\frac{\lambda}{\lambda -1}+\left(\frac{\lambda x}{\lambda -1}  -\frac{\lambda^2 }{\left(\lambda -1\right)^2}\right)\log \lambda,
	\\
	&&\mathfrak{D}_2\left(x;\lambda\right)=-\frac{2\lambda^3}{\left(\lambda-1\right)^2}+\frac{\lambda^2 \left(2x-1\right)}{\lambda-1}+\left(\frac{2\lambda^4}{\left(\lambda -1\right)^3}-\frac{2\lambda^3 x}{\left(\lambda -1\right)^2}+\frac{\lambda^2 x\left(x-1\right)}{\lambda-1}\right)
	\log \lambda,
\end{eqnarray*}
and so on.

\begin{figure}[H]
	\centering
	\includegraphics[width=0.6\textwidth]{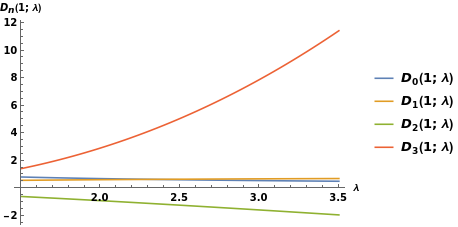} 
	\caption{\textmd{Plots of the polynomials $\mathfrak{D}_n\left(1;\lambda \right)$ for randomly selected special case when $\lambda\in\left[\frac{3}{2},\frac{7}{2}\right]$ with $n \in \{0,1,2,3\}$.}}
	\label{PlotsDPoly1a}
\end{figure}

By using the Implementation \ref{YlistDpoly} and the \texttt{Plot} command in Mathematica, we also give some two dimensional plots of the polynomials $\mathfrak{D}_n\left(x;\lambda \right)$ in Figure \ref{PlotsDPoly1a}-Figure \ref{PlotsDPoly3}.

\begin{figure}[H]
	\centering
	\includegraphics[width=0.6\textwidth]{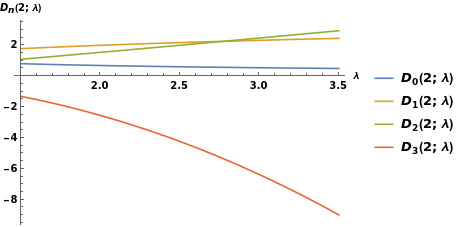} 
	\caption{\textmd{Plots of the polynomials $\mathfrak{D}_n\left(2;\lambda \right)$ for randomly selected special case when $\lambda\in\left[\frac{3}{2},\frac{7}{2}\right]$ with $n \in \{0,1,2,3\}$.}}
	\label{PlotsDPoly1b}
\end{figure}

\begin{figure}[H]
	\centering
	\includegraphics[width=0.6\textwidth]{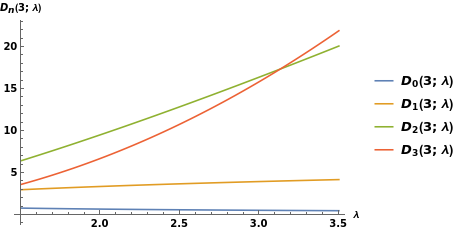} 
	\caption{\textmd{Plots of the polynomials $\mathfrak{D}_n\left(3;\lambda \right)$ for randomly selected case when $\lambda\in\left[\frac{3}{2},\frac{7}{2}\right]$ with $n \in \{0,1,2,3\}$.}}
	\label{PlotsDPoly1c}
\end{figure}

Figure \ref{PlotsDPoly1a}-Figure \ref{PlotsDPoly1c} shows the effects of the parameter $x$ on the graphs of the polynomials $\mathfrak{D}_n\left(x;\lambda \right)$.

Figure \ref{PlotsDPoly3} shows the effects of the parameter $\lambda$ on the graphs of the polynomials $\mathfrak{D}_n\left(x;\lambda \right)$.

\begin{figure}[H]
	\begin{subfigure}{\textwidth}
		\centering
		\includegraphics[width=0.65\textwidth]{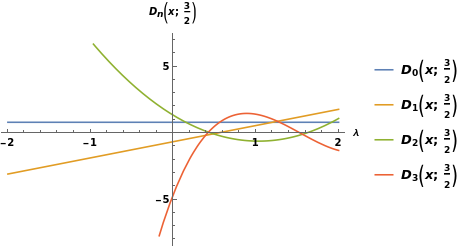} 
		\caption{}
	\end{subfigure}
	\\ 
	\begin{subfigure}{\textwidth}
		\centering
		\includegraphics[width=0.65\textwidth]{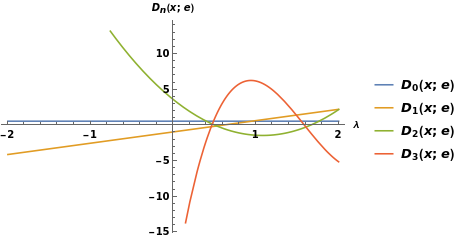}
		\caption{}
	\end{subfigure}
	\\
	\begin{subfigure}{\textwidth}
		\centering
		\includegraphics[width=0.65\textwidth]{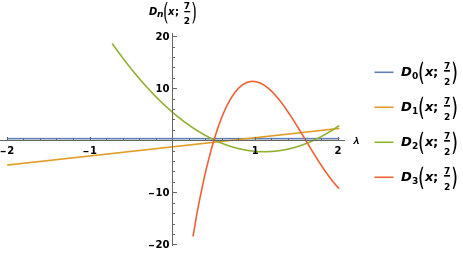} 
		\caption{}
	\end{subfigure}
	\caption{\textmd{Plots of the polynomials $\mathfrak{D}_n\left(x;\lambda \right)$ for the randomly selected special cases when $x\in\left[-2,2\right]$ and $n \in \{0,1,2,3\}$ \textbf{(a)} $\lambda=\frac{3}{2}$; \textbf{(b)} $\lambda=e$; \textbf{(c)} $\lambda=\frac{7}{2}$; \textbf{(d)} $\lambda=e^2$.}}
	\label{PlotsDPoly3}
\end{figure}

\begin{figure}[H]
	\ContinuedFloat
	\begin{subfigure}{\textwidth}
		\centering
		\includegraphics[width=0.65\textwidth]{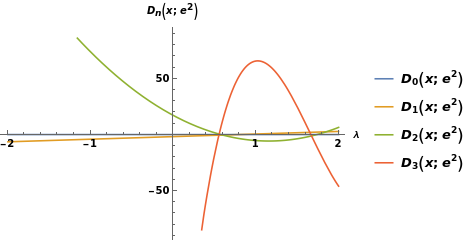} 
		\caption{}
	\end{subfigure}
	\caption{\textmd{Plots of the polynomials $\mathfrak{D}_n\left(x;\lambda \right)$ for the randomly selected special cases when $x\in\left[-2,2\right]$ and $n \in \{0,1,2,3\}$ \textbf{(a)} $\lambda=\frac{3}{2}$; \textbf{(b)} $\lambda=e$; \textbf{(c)} $\lambda=\frac{7}{2}$; \textbf{(d)} $\lambda=e^2$.}}
	\label{PlotsDPoly3}
\end{figure}

\section{Conclusion}
\label{Section-5}

In this paper, we present various identities and computation formulas containing not only the $\lambda$-Apostol-Daehee numbers and polynomials, but also Simsek numbers and polynomials, the Stirling numbers of the first kind, the Daehee numbers, and also the Chu-Vandermonde identity. Furthermore, by using functional equations containing the generating functions for the Cauchy numbers and the integrals of the generating functions for the $\lambda$-Apostol-Daehee numbers and polynomials, we also derive some identities and formulas for these numbers and polynomials. In addition, we give Mathematica implementation of a computation formula which computes the $\lambda$-Apostol-Daehee polynomials in terms of the polynomials $Y_n\left(x; \lambda\right)$. By the aid of the Mathematica implementation, we also give some plots which help the readers to analyze the behaviour of the $\lambda$-Apostol-Daehee polynomials for some randomly selected special cases of its parameters. As a conclusion, the results of this paper have the potential to affect many researchers conducting a research not only in computational mathematics, discrete mathematics and combinatorics, but also in other related fields.

For future studies, it is planned to investigate connections of the $\lambda$-Apostol-Daehee numbers with some special functions such as the Bernstein basis functions which possess many applications not only in approximation theory, but also in the construction of the Bezier curves widely used in computer-aided geometric design (\textit{cf}. \cite{AcikgozSerkan, ErkusDuman, Duman, Lorentz,FPTASimsek2013,SimsekBVP2013, SimsekHJMS2014, SimsekMMAS2015, Simsek Acikgoz} and also cited references therein).


\begin{thebibliography}{99}


\bibitem{AcikgozSerkan} Acikgoz, M., Araci, S., On generating function of the Bernstein polynomials, \textit{AIP Conf. Proc.}, 1281 (1) (2010), 1141, \url{https://dx.doi.org/10.1063/1.3497855}.

\bibitem{BonaBOOK} Bona, M., Introduction to Enumerative Combinatorics, The McGraw-Hill Companies Inc., New York, 2007.

\bibitem{CharalambidesBOOK} Charalambides, C. A., Enumerative Combinatorics, Chapman and Hall/ CRC Press Company, London, 2002.

\bibitem{ChoiFAR} Choi, J., Note on Apostol-Daehee polynomials and numbers, \textit{Far East J. Math. Sci.}, 101 (8) (2017), 1845--1857.

\bibitem{Comtet} Comtet, L., Advanced Combinatorics, D. Reidel Publication Company, Dordrecht-Holland/Boston-U.S.A., 1974.

\bibitem{Jordan} Jordan, C., Calculus of Finite Differences (2nd ed.), Chelsea Publishing Company, New York, 1950.

\bibitem{El-Desouky} El-Desouky, B. S., Mustafa, A., New results and matrix representation for Daehee and Bernoulli numbers and polynomials, \textit{Appl. Math. Sci. (Ruse)} 9 (2015), 3593--3610.

\bibitem{Gradimir} Cakic, N. P., Milovanovic, G. V., On generalized Stirling numbers and polynomials, \textit{Mathematica Balkanica} 18 (2004), 241--248.

\bibitem{Duman} Duman, O., Nonlinear Bernstein-type operators providing a better error estimation, \textit{Miskolc Math. Notes}, 15 (2) (2014), 393--400.

\bibitem{ErkusDuman} Erkus, E., Duman, O., Srivastava, H. M., Statistical approximation of certain positive linear operators constructed by means of the Chan-Chyan-Srivastava polynomials, \textit{Appl. Math. Comput.}, 182 (2006), 213--222.

\bibitem{KimDaehee} Kim, D. S., Kim, T., Daehee numbers and polynomials, \textit{Appl. Math. Sci. (Ruse)}, 7 (2013), 5969--5976.

\bibitem{KimEtAl2013} Kim, D. S., Kim, T., Lee, S.-H., Seo, J.-J., A Note on the lambda-Daehee Polynomials, \textit{Int. J. Math. Anal}, 7 (62) (2013), 3069--3080.

\bibitem{KimEtAl2016} Kim, T., Kim, D. S., Dolgy, D. V., Seo J.-J., Bernoulli polynomials of the second kind and their identities arising from umbral calculus, \textit{J. Nonlinear Sci. Appl.} 9 (2016), 860--869.

\bibitem{KucukogluAADM2019} Kucukoglu, I., Simsek, Y., On a family of special numbers and polynomials associated with Apostol-type numbers and poynomials and combinatorial numbers, \textit{Appl. Anal. Discrete Math.}, 13 (2019), 478--494.

\bibitem{KucukogluAxiom2019} Kucukoglu, I., Simsek, B., Simsek, Y., Generating Functions for New Families of Combinatorial Numbers and Polynomials: Approach to Poisson--Charlier Polynomials and Probability Distribution Function, \textit{Axioms}, 8 (4) (2019), 112, \url{https://dx.doi.org/10.3390/axioms8040112}.

\bibitem{Lorentz} Lorentz, G. G., Bernstein Polynomials, Chelsea Pub. Comp., New York, NY, USA, 1986.

\bibitem{MerliniCauchy} Merlini, D., Sprugnoli, R., Verri, M. C., The Cauchy numbers, \textit{Discrete Math.}, 306 (16) (2006), 1906--1920.

\bibitem{Qi} Qi, F., Explicit formulas for computing Bernoulli numbers of the second kind and Stirling numbers of the first kind, \textit{Filomat}, 
28 (2) (2014), 319--327.

\bibitem{ParkJCAA2016} Park, J.-W., On the $\lambda$-Daehee polynomials with $q$-parameter, \textit{J. Comput. Anal. Appl.}, 20 (1) (2016), 11--20.

\bibitem{Roman} Roman, S., The Umbral Calculus, Dover Publ. Inc., New York, 2005.

\bibitem{FPTASimsek2013} Simsek, Y., Functional equations from generating functions: A novel approach to deriving identities for the Bernstein basis functions, \textit{Fixed Point Theory Appl.}, Article number: 80 (2013), 1--13.

\bibitem{SimsekFPTA} Simsek, Y., Generating functions for generalized
Stirling type numbers, array type polynomials, Eulerian type polynomials and
their alications, \textit{Fixed Point Theory Appl.}, 87 (2013), 343--355.

\bibitem{SimsekBVP2013} Simsek, Y., Unification of the Bernstein-type polynomialsand their applications, \textit{Bound. Value Probl.}, Article number: 56 (2013), 1--15.

\bibitem{SimsekHJMS2014} Simsek, Y., Generating functions for the Bernstein type polynomials: A new approach to deriving identities and applications for the polynomials, \textit{Hacet. J. Math. Stat.}, 43 (1) (2014), 1--14.

\bibitem{SimsekMMAS2015} Simsek, Y., Analysis of the Bernstein basis functions: an approach to combinatorial sums involving binomial coefficients
and Catalan numbers, \textit{Math. Methods Appl. Sci.}, 38 (14) (2015), 3007--3021.

\bibitem{simsekCogent} Simsek, Y., Analysis of the $p$-adic $q$-Volkenborn integrals: An approach to generalized Apostol-type special numbers and polynomials and their applications, \textit{Cogent Math. Stat.} 1269393 (2016), \url{https://dx.doi.org/10.1080/23311835.2016.1269393}.

\bibitem{SimsekASCM2016} Simsek, Y., Apostol type Daehee numbers and polynomials, \textit{Adv. Stud. Contemp. Math. (Kyungshang)}, 26 (3) (2016), 555--566.

\bibitem{SimsekASCM2017} Simsek, Y., Identities on the Changhee numbers and Apostol-type Daehee polynomials, \textit{Adv. Stud. Contemp. Math. (Kyungshang)}, 27 (2) (2017), 199--212.

\bibitem{simsekTJM} Simsek, Y., Construction of some new families of Apostol-type numbers and polynomials via Dirichlet character and $p$-adic $q$-integrals, \textit{Turk. J. Math.}, 42 (2018), 557--577.

\bibitem{ysimsekMTJPAM} Simsek, Y., Explicit formulas for $p$-adic integrals: Approach to $p$-adic distributions and some families of special	numbers and polynomials, \textit{Montes Taurus J. Pure Appl. Math.}, 1 (1) (2019), 1--76.

\bibitem{Simsek Acikgoz} Simsek, Y., Acikgoz, M., A new generating function of ($q$-) Bernstein-type polynomials and their interpolation
function, \textit{Abstr. Appl. Anal.}, 769095 (2010), \url{https://dx.doi.org/10.1155/2010/769095}.

\bibitem{SimsekYardimci2016} Simsek, Y., Yardimci, A., Applications on the Apostol-Daehee numbers and polynomials associated with special numbers, polynomials, and $p$-adic integrals, \textit{Adv. Difference Equ.}, 308 (2016), \url{https://dx.doi.org/10.1186/s13662-016-1041-x}.

\bibitem{SimsekKucukoglu2020Chapter} Simsek, Y., Kucukoglu, I., Some Certain Classes of Combinatorial Numbers and Polynomials Attached to Dirichlet Characters: Their Construction by $p$-adic Integration and Applications to Probability Distribution Functions, A Chapter in the Book: Mathematical Analysis in Interdisciplinary Research, Th. M. Rassias, I. N. Parasidis and E. Providas (Eds.), Springer International Publishing, Springer Nature Switzerland AG, In press.

\bibitem{KucukogluJNT2017} Srivastava, H. M., Kucukoglu, I., Simsek, Y., Partial differential equations for a new family of numbers and polynomials unifying the Apostol-type numbers and the Apostol-type polynomials, \textit{J. Number Theory}, 181 (2017), 117--146.

\bibitem{WolframCloud} Wolfram Research Inc., Mathematica Online (Wolfram Cloud), Champaign,
IL, 2020, \url{https://www.wolframcloud.com}

\end{thebibliography}
\end{document}